\documentclass{article}
\title{\bf What an infra-nilmanifold endomorphism really should be \ldots}
\author{Karel Dekimpe\\
K.U.Leuven Campus Kortrijk, Universitaire Campus, B-8500 Kortrijk.}
\long\def\symbolfootnote[#1]#2{\begingroup%
\def\thefootnote{\fnsymbol{footnote}}\footnote[#1]{#2}\endgroup} 
\usepackage{amssymb}
\usepackage{a4wide}
\usepackage{xypic}
\newcommand{\Z}{{\mathbb Z}}
\newcommand{\Q}{{\mathbb Q}}
\newcommand{\R}{{\mathbb R}}

\newcommand{\Aut}{{\rm Aut}}
\newcommand{\Aff}{{\rm Aff}}
\newcommand{\Diff}{{\rm Diff}}
\newcommand{\Isom}{{\rm Isom}}
\newcommand{\GL}{{\rm GL}}
\newcommand{\semi}{{\rtimes}}
\newcommand{\qed}{\hfill {$\Box$}}
\newcommand{\lien}{{\mathfrak n}} 
\newtheorem{theorem}{Theorem}[section]
\newtheorem{proposition}[theorem]{Proposition}
\newtheorem{corollary}[theorem]{Corollary}
\newtheorem{definition}[theorem]{Definition}
\newtheorem{lemma}[theorem]{Lemma}
\newtheorem{remark}[theorem]{Remark}
\newtheorem{question}[theorem]{Open Question}
\begin{document}
\maketitle
\begin{abstract}
Infra-nilmanifold endomorphisms were introduced in the late sixties. They 
play a very crucial role in dynamics, especially when studying expanding 
maps and Anosov diffeomorphisms.
However, in this note we will explain that the two main results in this area 
are based on a false result and that although we can repair one of these two theorems, 
there remains doubt on the correctness of the other one.
Moreover, we will also show that the notion of an 
infra-nilmanifold endomorphism itself has not always been interpreted in the same way. \\
Finally, we define a slightly more general concept of the notion of an 
infra-nilmanifold endomorphism and explain why this is really the right concept to work with.%
\symbolfootnote[0]{{\bf Keywords:} Anosov diffeomorphism, Expanding map, Nilpotent Lie group, Infra-nilmanifold, Affine map}
\symbolfootnote[0]{{\bf MSC 2010:} 37D20, 20F34}
\end{abstract}
\section{Introduction}

The notion of an infra-nilmanifold endomorphism appears for the first time 
in the proceedings of the symposium in pure mathematics of the American Mathematical
Society held in 1968 in Berkeley (\cite{fran70-1,hirs70-1,shub70-1}). 

Nowadays, when using the term infra-nilmanifold endomorphism, most people 
refer to the paper of J.~Franks (\cite{fran70-1}), although J.~Franks himself 
in that same paper (and also M.~Shub in \cite{shub70-1}), attributes this terminology to 
M.W.~Hirsch (\cite{hirs70-1}). However, it is immediately clear and this will 
also be explained in the next section, that the definition used
by J.~Franks and M.~Shub is not equivalent to the one used by M.W.~Hirsch. 

The results of (\cite{fran70-1}) have played an important role in dynamics, especially 
in the study of Anosov diffeomorphisms and expanding maps. 
On the one hand, it was a crucial ingredient in the result of M.~Gromov (See geometric corollary on page 
55 of \cite{grom81-1}) stating that an expanding map on an arbitrary compact manifold is topologically conjugate 
to an infra-nilmanifold endomorphism.

On the other hand, A.~Manning (\cite{mann74-1}), using  the concept of an infra-nilmanifold endomorphism as introduced in (\cite{fran70-1}),
showed that any Anosov diffeomorphism of an infra-nilmanifold $M$ is 
topologically conjugate to a so called hyperbolic infra-nilmanifold automorphism. 

\medskip

Unfortunately, some results of \cite{fran70-1} depend on a ``theorem'' of L.~Auslander 
(\cite[Theorem 2]{ausl60-1}) which is not correct (not only the proof, but the 
statement of Auslander's  theorem is wrong). We will explain this in the next section.
Moreover, although most of the arguments in \cite{fran70-1} which are based on Auslander's wrong result can
be restored using a modified version of it (see Corollary~\ref{correctAusl} below), there remain some subtle
problems with the definition of the concept of an infra-nilmanifold endomorphism as given in (\cite{fran70-1}).

\medskip

The aim of this note is to show that, even with a correct definition of an infra-nilmanifold endomorphism, 
both the proofs of the result of A.~Manning and of M.~Gromov are not correct, because they are heavily based 
on a wrong result (\cite[Proposition 3.5]{fran70-1}) of the paper of Franks. 
As both of these results are often referred to, I will point out as detailed as possible, where the
problems in the work of L.~Auslander and of J.~Franks are situated and
how this has it's influence in the work of M.~Gromov and A.~Manning. Moreover, I will give an example of an expanding 
map and of an Anosov diffeomorphism on a given infra-nilmanifold which are not topologically conjugate to an infra-nilmanifold
endomorphism of that infra-nilmanifold. Fortunately, by the work of K.B.~Lee and F.~Raymond (\cite{lr85-1}), who were,
up to my knowledge, 
the first to discover the problems in the work of  L.~Auslander, it is rather easy to define 
a slightly broader concept of the notion of infra-nilmanifold endomorphism, namely the class of affine endomorphisms,
which is more suited to study self maps of infra-nilmanifolds. 
We will show that using this broader concept the result of M.~Gromov on expanding maps can be repaired, but one has to be very careful
with the precise interpretation of the statement.
On the other hand, although it is also to be expected that A.~Manning's result might be repaired, I haven't
been able to prove this in its full generality yet.

\section{Infra-nilmanifolds and endomorphisms of their fun\-da\-men\-tal groups}

Let $N$ be a connected and simply connected nilpotent Lie group and let $\Aut(N)$ 
be the group of continuous automorphisms of $N$. Then $\Aff(N)=N\semi \Aut(N)$ acts on 
$N$ in the following way:
\[ \forall (n,\alpha) \in \Aff(N),\;\forall x \in N: (n,\alpha)\cdot x = n \alpha(x).\]
So an element of $\Aff(N)$ consists of a translational part $n\in N$ and a linear part $\alpha\in \Aut(N)$
(as a set $\Aff(N)$ is just $N\times \Aut(N)$) and $\Aff(N)$ acts on $N$ by first applying the 
linear part and then multiplying on the left by the translational part). In this way, $\Aff(N)$ can 
also be seen as a subgroup of $\Diff(N)$. 

Now, let $C$ be a compact subgroup of $\Aut(N)$ and consider any torsion free discrete 
subgroup $\Gamma$ of $N\semi C$, such that the orbit space $\Gamma\backslash N$ is compact.
Note that $\Gamma$ acts on $N$ as being also a subgroup of $\Aff(N)$. 

The action of $\Gamma$ on $N$ will be free and properly discontinuous, so $ \Gamma\backslash N$
is a manifold, which is called an infra-nilmanifold. 
It follows from the (correct) Theorem 1 of L.~Auslander in (\cite{ausl60-1}), that 
$\Gamma \cap N$ is a uniform lattice of $N$ and that $\Gamma/(\Gamma\cap N)$ is a finite 
group. This shows that the fundamental group of an infra-nilmanifold   $ \Gamma\backslash N$
is virtually nilpotent (i.e.\ has a nilpotent normal subgroup of finite index). In fact 
$\Gamma \cap N$ is a maximal nilpotent subgroup of $\Gamma$ and it is the only normal 
subgroup of $\Gamma$ with this property. (This also follows from \cite{ausl60-1}).

If we denote by $p:N\semi C \rightarrow C$ the natural projection on the 
second factor, then $p(\Gamma)\cong  \Gamma/(\Gamma\cap N)$. Let $F$ denote this finite 
group $p(\Gamma)$, then we will refer to $F$ as being the holonomy group of $\Gamma$ (or
of the infra-nilmanifold $\Gamma\backslash N$). It follows that 
$\Gamma \subseteq N\semi F$.  In case $F=1$, so $\Gamma\subseteq N$, the manifold 
$N\backslash G$ is a nilmanifold. Hence, any infra-nilmanifold $\Gamma\backslash N$ 
is finitely covered by a nilmanifold $(\Gamma\cap N)\backslash N$. This also explains
the prefix ``infra''.

When the Lie group $N$ is abelian, so $N$ is the additive group $\R^n$ for some $n$,
it is enough to consider the case $C=O(n)$, the orthogonal group, because $O(n)$ is a maximal 
compact subgroup of $\Aut(\R^n)=\GL_n(\R)$ and so any other compact subgroup is 
conjugate to a subgroup of $O(n)$. It follows that in this situation $N\semi C = \R^n\semi O(n)$ 
is the group of isometries of Euclidean space $\R^n$. In this setting, the infra-nilmanifolds
are compact flat Riemannian manifolds and the nilmanifolds are just tori.

\begin{remark}
Many authors (E.g.~see \cite{fran70-1}, \cite{hirs70-1}) 
start from discrete subgroups of $N\semi F$ for various finite groups $F$ to 
define the notion of an infra-nilmanifold. The discussion above shows that this is not a 
restriction.\\
In (\cite{grom81-1}) and (\cite{hira88-1}), an infra-nilmanifold is defined as a quotient $\Gamma\backslash N$, where
$\Gamma$ is a subgroup of the whole affine group $\Aff(N)$ acting freely and properly discontinuously 
on $N$. This is not a correct definition, for in this case, the linear parts do not have to form 
a finite group and hence $\Gamma$ need not be a virtually nilpotent group. 
As an example, let $\varphi:\Z \rightarrow \Aut(\Z^2)$ be any morphism and regard 
$\varphi(z)$ as being a $2\times 2$--matrix. Then, 
\[ \Gamma= \left\{ \left( \left(\begin{array}{c} x \\ y \\ z
\end{array}\right), \left( \begin{array}{cc} \varphi(z) & 0 \\ 0 & 1\end{array}\right) \right)\;|\;
x,y,z\in \Z \right\} \]
is a subgroup of $\Aff(\R^3)$ acting freely and properly discontinuously on $\R^3$. The group
$\Gamma$ is isomorphic to the semi-direct product group $\Z^2 \semi \Z$, where
the action of $\Z$ on $\Z^2$ is given via $\varphi$. Such a group is often not virtually 
nilpotent. E.g.\ there is a unique morphism $\varphi:\Z \rightarrow \Aut(\Z^2)$, with 
$\varphi(1)=\left( \begin{array}{cc} 2 & 1 \\ 1 & 1\end{array}\right)$. The corresponding 
group $\Z^2\semi \Z$ is not virtually nilpotent. Actually, the manifolds which are obtained
in this way are called complete affinely flat manifolds (see \cite{miln77-1}). 
\end{remark} 

Let us now discuss why Theorem 2 of \cite{ausl60-1} is not correct. In fact, L.~Auslander proves
this theorem as a generalization of the second Bieberbach theorem. Unfortunately, even L.~Auslander's formulation
of this second Bieberbach theorem is not correct. This was first observed, without further explanation, by K.B.~Lee and 
F.~Raymond in \cite{lr85-1}. As this theorem plays an important role in the work of J.~Franks, I will
explain in full detail what goes wrong and what can be saved. 

We recall the statement of Auslander's theorem using the notations we introduced above.
\begin{quote}{\bf Formulation of Theorem 2 in \cite{ausl60-1}}\\
Let $\Gamma_1$ and $\Gamma_2$ be discrete uniform subgroups of $N\semi C$. Let $\psi:\Gamma_1
\rightarrow \Gamma_2$ be an isomorphism. Then $\psi$ can be uniquely extended to a continuous  automorphism 
$\psi^\ast$ of $N\semi C$ onto itself.
\end{quote}

It is very easy to produce a counterexample to this statement. In fact, the statement is almost never 
correct. Let $N=\R^2$ the additive group and   
$C=O(2)$. Let $\Gamma_1=\Gamma_2=\Z^2$ and let $\psi\in \Aut(\Z^2)$ be the automorphism represented by the 
matrix $A=\left( \begin{array}{cc} 2 & 1 \\ 1 & 1\end{array}\right)$ (almost any matrix will do). 
Now assume that $\psi$ extends to a continuous automorphism $\psi^\ast$ of $\R^2 \semi O(2)$. The group $\R^2$ (seen 
as a subgroup of $\R^2\semi O(2)$) is 
normal and maximal abelian and is the unique subgroup of $\R^2\semi O(2)$ satisfying this condition, so 
we must have that $\psi^\ast(\R^2)= \R^2$. It follows that the restriction of $\psi^\ast$ to $\R^2$ is 
the linear map, given by the matrix $A$.  So $\psi^\ast(r,1)=(Ar,1)$ for all $r\in \R^2$. (Here $1$ denotes the 
trivial automorphism of $\R^2$ or the $2\times 2$ identity matrix)

Now let $B\in O(2)$, so $(0,B)\in \R^2\semi O(2)$, and assume that $\psi^\ast(0,B)=(b,B')$ for 
some $b\in \R^2$ and some $B'\in O(2)$. Let us perform a small computation, where $r\in \R^2$ is arbitrary:
\begin{eqnarray*}
\psi^\ast((0,B)(r,1)(0,B)^{-1}) & = & \psi^\ast(0,B) \psi^\ast(r,1) \psi^\ast(0,B)^{-1}\\
& \Downarrow & \\
\psi^\ast(Br,1) & = & (b,B') ( Ar ,1) (-B'^{-1}b, B'^{-1})\\
& \Downarrow & \\
(ABr,1) & = & (B'A r,1)\\
\end{eqnarray*}
As this holds for any $r$ we must have that $AB=B'A$, or $B'=ABA^{-1}$. It is now trivial to see 
such a $B'$ does not have to belong to $O(2)$. E.g.~when $B=\left(\begin{array}{cc} -1 & 0 \\0 & 1 \end{array}\right)$.
We have that $B'=  \left(\begin{array}{cc} -3 & 4 \\-2 & 3 \end{array}\right)\not \in O(2)$. 
We can conclude that $\psi$ does not extend to a continuous morphism of $\R^2\times O(2)$, contradicting the 
statement made by L.~Auslander. At this point I want to remark that the ``proof'' of Auslander is very short
and does not make any sense to me, so it is difficult to point out where exactly the error is situated in his argument.

\medskip

A correct formulation of a generalization of the second Bieberbach theorem is given 
in \cite{lr85-1}.
\begin{theorem}\label{correctionAuslander} (\cite{lr85-1}, see also \cite[page 16]{deki96-1})
Let $N$ be a connected and simply connected nilpotent Lie group and $C$ a compact 
subgroup of $\Aut(N)$.
Let $\Gamma_1$ and $\Gamma_2$ be two discrete and uniform subgroups of $N\semi C$ and 
assume that $\psi:\Gamma_1 \rightarrow \Gamma_2$ is an isomorphism, then 
there exists a $\alpha\in \Aff(N)$ such that
\[ \forall \gamma\in \Gamma_1: \; \psi(\gamma)= \alpha \gamma \alpha^{-1}.\]
\end{theorem}
So, any isomorphism between the groups $\Gamma_1$ and $\Gamma_2$ is induced by a conjugation inside $\Aff(N)$.

\medskip

At this point, I would like to mention a corollary, which can be seen as a fix to the 
false statement of L.~Auslander.

\begin{corollary}\label{correctAusl}
Let $N$ be a connected and simply connected nilpotent Lie group and $C$ a compact 
subgroup of $\Aut(N)$ and let $\Gamma$ be a discrete and uniform subgroup of $N\semi C$.
Let $p:N\semi C\rightarrow C$ denote the natural projection. \\
If $\psi:\Gamma\rightarrow \Gamma$ is a monomorphism, then $p(\Gamma)=p(\psi(\Gamma))$. Moreover, in this
case $\psi$ extends to an automorphism $\psi^\ast$ of  $N\semi p(\Gamma)$, such that $\psi^\ast(N)=N$.
\end{corollary}

\underline{Proof:} 
$\psi$ is an isomorphism from $\Gamma$ onto $\psi(\Gamma)$, so by Theorem~\ref{correctionAuslander},
$\psi$ can be realized as a conjugation, say by $\alpha\in \Aff(N)$, inside $\Aff(N)$.
As $N$ is a normal subgroup of $\Aff(N)$, we have that $\alpha N \alpha^{-1} = N $. On the other 
hand, we also have that $\alpha \Gamma \alpha^{-1} \subseteq \Gamma$. Therefore,
$\alpha (N \semi p(\Gamma))\alpha^{-1} = \alpha N\Gamma \alpha^{-1} \subseteq N\Gamma = N\semi p(\Gamma)$.
In fact, we can see that $\alpha (N\semi p(\Gamma)) \alpha^{-1}= N\semi p(\Gamma)$ (and not a 
proper subset of it). To prove this, we must show that for any $\mu\in p(\Gamma)$, there is 
a $n \in N$, with $(n,\mu)\in \psi(\Gamma)= \alpha \Gamma \alpha^{-1}$.
This is however easy, because any morphism $\psi$ of $\Gamma$ induces a morphism 
\[ \bar\psi: p(\Gamma)=\Gamma/(\Gamma\cap N) \rightarrow p(\Gamma)= \Gamma/(\Gamma\cap N).\]
Now, as $\psi$ is conjugation with an element $\alpha\in \Aff(N)$, it is easy to see 
that $\bar\psi$ is conjugation with the linear part of $\alpha$ in $\Aut(N)$. Therefore,
$\bar\psi$ is bijective, 
showing that $p(\psi(\Gamma))= p(\Gamma)$ and $\alpha (N\semi p(\Gamma)) \alpha^{-1}= N\semi p(\Gamma)$.
The proof now finishes by taking $\psi^*$ to be conjugation with $\alpha$ inside $\Aff(N)$. \qed

\section{Infra-nilmanifold endomorphisms}
In this section, we will discuss the notion of an infra-nilmanifold endomorphism
as introduced by M.W.~Hirsch (\cite{hirs70-1}) and by J.~Franks (\cite{fran70-1}).

To do this, we fix an infra-nilmanifold $\Gamma\backslash N$, so $N$ is a connected and 
simply connected nilpotent Lie group and $\Gamma$ is a torsion free, uniform discrete 
subgroup of $N\semi F$, where $F$ is a finite subgroup of $\Aut(N)$. We will assume
that $F$ is the holonomy group of $\Gamma$ (so for any $\mu\in F$, there exists a $n\in N$ such 
that $(n,\mu)\in \Gamma$).  

In what follows, we will identify $N$ with the subgroup $N\times \{1\}$ of $N\semi \Aut(N)=\Aff(N)$,
$F$ with the subgroup $\{1\} \times F$ and $\Aut(N)$ with the subgroup $\{1\}\times \Aut(N)$.
Hence, we can say that an element of $\Gamma$ is of the form $n\mu$ for some $n\in N$ and some 
$\mu\in F$. Also, any element of $\Aff(N)$ can uniquely be written as a product 
$n\psi$, where $n\in N$ and $\psi\in \Aut(N)$. The product in $\Aff(N)$ is then 
given as 
\[ \forall n_1,n_2\in N,\;\forall \psi_1,\psi_2\in \Aut(N):\; n_1 \psi_1 n_2 \psi_2= 
n_1 \psi_1(n_2) \psi_1\psi_2.\]

\medskip

We will first look at the way M.W.~Hirsch introduced the notion of an infra-nilmanifold endomorphism. 
Actually, Hirsch defines endomorphisms on a larger class of spaces, called infra homogeneous spaces, but we 
immediately specialise to the case of infra-nilmanifolds.  

M.W.~Hirsch starts with a given automorphism $\varphi$ of the Lie group $N\semi F$, 
with $\varphi(F)=F$. Note that we also have that  $\varphi(N)=N$, because $N$ is the
connected component of the identity element in $N\semi F$.  Before we continue, let us 
give a description of these automorphisms. 

\begin{lemma}\label{DucDeBrabant} Let $N$ be a connected, simply connected nilpotent Lie group and $F$ be a  finite 
subgroup of $\Aut(N)$. Let $\varphi$ be an automorphism of $N\semi F$ with $\varphi(F)=F$ and denote 
by $\psi\in \Aut(N) $ the restriction of $\varphi$ to $N$, then
\[\forall x\in N\semi F:\; \varphi (x) = \psi x \psi^{-1} \]
where $\psi x \psi^{-1}$ is  a conjugation in the group $\Aff(N)$. 
\end{lemma}
\underline{Proof:} 
For any $n\in N$ and any $\xi\in \Aut(N)$, the equality 
$\xi(n)= \xi n \xi^{-1}$ is valid, where $\xi n \xi^{-1}$ is a conjugation in $\Aff(N)$. 
So, we also have that 
\begin{equation}\label{step1}
\varphi(n)= \psi(n)= \psi n\psi^{-1}.
\end{equation}
Let us now consider an element $\mu\in F$. 
For any $n\in N$, we have the following equation in the group $N\semi F$:
\[ \mu(n) = \mu n \mu^{-1}.\]
By applying $\varphi$ to both sides of this equation, we find that 
\begin{eqnarray*}
\varphi(\mu(n)) & = & \varphi(\mu) \psi(n) \varphi(\mu)^{-1}\\
& \Downarrow & \\
\psi(\mu(n)) & = & \varphi(\mu) (\psi(n))
\end{eqnarray*}
Since this holds for any $n\in N$, we have that 
\[ \psi \circ \mu = \varphi(\mu) \circ \psi\] 
showing that 
\begin{equation}\label{step2}
\varphi(\mu) = \psi \mu \psi^{-1}.
\end{equation}
Now, combining (\ref{step1}) and (\ref{step2}) we find that for $x=n\mu$, with $n\in N$ and 
$\mu\in F$:
\[ \varphi(x) = \varphi(n)\varphi(\mu)= \psi n \psi^{-1} \psi \mu \psi^{-1}= \psi x \psi^{-1},\] 
which finishes the proof. \qed

\medskip

Now, let $\varphi$ still be  an automorphism of $N\semi F$ with $\varphi(F)=F$ and assume that 
$\varphi(\Gamma) \subseteq \Gamma$, where $\Gamma$ is a torsion free, discrete and uniform 
subgroup of $N\semi F$.  Now, let $\gamma=m\mu$ be any element of $\Gamma$, where 
$m\in N$ and $\mu \in F$. We denote the action of $\Gamma$ on $n\in N$ by 
$\gamma\cdot n$, so $\gamma \cdot n = m \mu(n)$.
Now we compute that 
\begin{eqnarray*}
\varphi(\gamma\cdot n) & = & \varphi(m \mu(n)) \\
& = & \varphi(m) \varphi(\mu(n) ) \\
& = & \varphi(m) \varphi(\mu)(\varphi(n))  \\
& = & \varphi(m\mu)\cdot \varphi(n)\\
& = & \varphi(\gamma) \cdot \varphi(n).
\end{eqnarray*}

We are now ready to introduce the notion of an infra-nilmanifold endomorphism.
\begin{definition}[Infra-nilmanifold endomorphism following Hirsch]\label{defHirsch}
Let $N$ be a connected and simply connected nilpotent Lie group and $F\subseteq \Aut(N)$ a finite 
group. Assume that $\Gamma$ is a torsion free, discrete and uniform subgroup of $N\semi F$.
Let $\varphi:N\semi F\rightarrow N\semi F$ be an automorphism, such that $\varphi(F)=F$ and 
$\varphi(\Gamma)\subseteq \Gamma$, then, the map  
\[ \bar\varphi:\; \Gamma\backslash N \rightarrow \Gamma\backslash N :\;
\Gamma\cdot  n \mapsto \Gamma\cdot \varphi(n),\]
is the infra-nilmanifold endomorphism induced by $\varphi$. In case $\varphi(\Gamma)=\Gamma$, we 
call $\bar\varphi$ an infra-nilmanifold automorphism.
\end{definition}
In the definition above, $\Gamma\cdot  n$ denotes the orbit of $n$ under the action of $\Gamma$.
The computation above shows that $\bar\varphi$ is well defined. Note that 
infra-nilmanifold automorphisms are diffeomorphisms, while in general an infra-nilmanifold endomorphism 
is a self-covering map. 

\begin{remark}
It is easy to check that for an infra-nilmanifold endomorphism $\bar\varphi$, the induced 
morphism $\bar\varphi_\sharp$ on the fundamental group $\Pi_1(\Gamma\backslash N, x_0)\cong \Gamma$ is exactly
the restriction of $\varphi$ to $\Gamma$ (see also Proposition~\ref{laternodig} below and note that
one can always choose as basepoint $x_0=\Gamma\cdot e$, the orbit of the identity element of $N$). 
By Lemma~\ref{DucDeBrabant} we know that $\bar\varphi_\sharp$ is induced
by a conjugation with an automorphism inside $\Aff(N)$. On the other hand, Theorem~\ref{correctionAuslander} shows
that in general an injective endomorphism of $\Gamma$ is induced by a conjugation with a 
general element of $\Aff(N)$ and not necessarily by an automorphism. This already indicates that 
there might exist (interesting) diffeomorphisms and self--covering maps of an infra-nilmanifold
which are not even homotopic to an infra-nilmanifold endomorphism. Further on, we will 
explicitly construct such examples and obtain an Anosov diffeomorphism (resp.\ an expanding map) of
an infra-nilmanifold
which is not homotopic to an infra-nilmanifold automorphism (resp.~infra-nilmanifold endomorphism) of that 
infra-nilmanifold.
\end{remark}

As already indicated above, we will also consider the definition of an infra-nilmanifold endomorphism 
as introduced by J.~Franks in \cite[page 63]{fran70-1}, the definition which is in fact most often referred to. 
Using our notation introduced above, 
J.~Franks writes that when $\varphi:N\semi F\rightarrow N\semi F$ is an automorphism for which 
$\varphi(\Gamma)\subseteq\Gamma$ and $\varphi(N)=N$, it induces a map
\[ \bar\varphi: \Gamma\backslash N \rightarrow \Gamma \backslash N.\]
(In fact, J.~Franks requires that $\varphi(\Gamma)=\Gamma$ and not that it is only a subgroup, but I 
believe this is a typo).\\
It is this kind of maps that he calls infra-nilmanifold endomorphisms.
As J.~Franks does not impose the condition that $\varphi(F)=F$, this seems to be a generalization 
of the  notion introduced by M.W.~Hirsch. 
Exactly the same definition was given by M.~Shub in \cite[page 274]{shub70-1} (without the typo).  

\medskip

Unfortunately, there seems to be a problem with this definition. It is {\bf not true} that 
the map $\bar\varphi: \Gamma\backslash N \rightarrow \Gamma \backslash N: \Gamma \cdot n \mapsto
\Gamma \cdot \varphi(n)$ is in general well defined. As many authors refer to the work of J.~Franks when talking 
about infra-nilmanifold endomorphisms, we give a detailed example to show where it goes wrong.

Let $N=\R^3$, the additive group. We let $F\cong \Z_2\oplus\Z_2\subseteq \GL_3(\R)$ be the group 
with elements
\[ 1=\left(\begin{array}{ccc} 1 & 0 & 0 \\ 0 & 1 & 0 \\ 0 & 0 & 1 \end{array} \right),\;
\alpha =\left(\begin{array}{ccc} 1 & 0 & 0 \\ 0 & -1 & 0 \\ 0 & 0 & -1 \end{array} \right),
\beta=\left(\begin{array}{ccc} -1 & 0 & 0 \\ 0 & 1 & 0 \\ 0 & 0 & -1 \end{array} \right),\;
\alpha\beta =\left(\begin{array}{ccc} -1 & 0 & 0 \\ 0 & -1 & 0 \\ 0 & 0 & 1 \end{array} \right).\]
Moreover, we pick 
\[ a=\left(\begin{array}{c} \frac{1}{2} \\ 0 \\ 0 \end{array} \right) \mbox{ and }
b = \left(\begin{array}{c} 0 \\ \frac{1}{2} \\ \frac12 \end{array} \right).\]

Let $A=(a,\alpha) \in \Aff(\R^3)$ and $B=(b,\beta)\in \Aff(\R^3)$ and consider 
the group $\Gamma\subseteq \R^3\semi F$ to be the group generated by $\Z^3\cup \{ A,B\}$.
Then $\Gamma$ is a torsion free, uniform discrete subgroup of $\R^3\semi F$. In fact 
$\Gamma\backslash \R^3$ is the well known Hantsche-Wendt manifold with fundamental 
group $\Gamma$.

Now, let $\delta=\left( \begin{array}{ccc} 1 & 0 & 0 \\ 0 & 0 & 1\\ 0 & 1 & 0 \end{array}\right)$ and 
$\displaystyle d=\left( \begin{array}{c} \frac{1}{4} \\ 0 \\ 0 \end{array} \right) $ and 
take $D=(d,\delta) \in \Aff(\R^3)$.
Let $\varphi: \Aff(\R^3) \rightarrow \Aff(\R^3): X \mapsto DXD^{-1}$ be the inner automorphism 
determined by $D$. A calculation shows that $\varphi (\R^3\semi F) =\R^3 \semi F$, so 
$\varphi$ restricts to an automorphism $\R^3\semi F$ for which of course $\varphi(\R^3)=\R^3$ 
(but $\varphi(F)\neq F$!). Moreover,
\[ \varphi(\Z^3)=\Z^3\subseteq \Gamma,\; \varphi(A) = A \in \Gamma \mbox{ and } 
\varphi(B)= \left( \begin{array}{c} 0 \\ 1 \\ 1\end{array} \right) A B \in \Z^3 \Gamma=\Gamma. \]
So $\varphi(\Gamma) \subseteq \Gamma$ (in fact equality holds). I claim that in this case
the map $\bar\varphi$ is not well defined. To prove this claim, we need to provide a $n\in \R^3$ and
a $\gamma\in \Gamma$, such that 
\[\Gamma\cdot \varphi(n) \neq \Gamma \cdot \varphi(\gamma \cdot n).\]
Let
\[ n=   \left( \begin{array}{c} \frac13 \\ \frac13 \\ \frac13 \end{array} \right) \mbox{ and }
\gamma= B, \mbox{ then } \gamma\cdot n=  \left( \begin{array}{c} -\frac13 \\ \frac56 \\ \frac16 \end{array} \right).\]
It follows that 
\[ \varphi(n)=\left( \begin{array}{c} \frac13 \\ \frac13 \\ \frac13 \end{array} \right)
\mbox{ and }
\varphi(\gamma\cdot n)= \left( \begin{array}{c} -\frac13 \\ \frac16 \\ \frac56 \end{array} \right).\]
To check that $\Gamma \cdot \varphi(\gamma\cdot n) \neq \Gamma \cdot \varphi(n)$, it suffices to check that 
$\varphi(\gamma\cdot  n) \not \in \Gamma \cdot \varphi(n)$, or that $\varphi(\gamma\cdot n) \neq \gamma'\cdot \varphi(n)$ for 
any $\gamma'\in \Gamma$. Any $\gamma'\in\Gamma$ can uniquely be written in one of the following ways:
\[ \gamma'_1 = z,\; \gamma'_2 = z A,\; \gamma'_3 = z B\mbox{ or } 
\gamma'_3= z AB,\mbox{ with } z= \left( \begin{array}{c} z_1 \\ z_2 \\ z_3 \end{array} \right)\in \Z^3.\]
Computing $\gamma'\cdot \varphi(n)$ in each of these four case, we obtain:
\[ \gamma'_1 \cdot \varphi(n) = \left( \begin{array}{c} z_1 + \frac13 \\ z_2 + \frac13 \\ z_3 +\frac13 \end{array} \right)\;
   \gamma'_2 \cdot \varphi(n) = \left( \begin{array}{c} z_1 + \frac56 \\ z_2 - \frac13 \\ z_3 -\frac13 \end{array} \right)\;
   \gamma'_3 \cdot \varphi(n) = \left( \begin{array}{c} z_1 - \frac13 \\ z_2 + \frac56 \\ z_3 +\frac16 \end{array} \right)\;
   \gamma'_4 \cdot \varphi(n) = \left( \begin{array}{c} z_1 + \frac16 \\ z_2 - \frac56 \\ z_3 -\frac16 \end{array} \right).
\]
It is obvious that none of these expressions equals $\varphi(\gamma\cdot n)$, proving the claim. \qed

\medskip
 
At the end of this section, we want to explain that in a certain sense, the definition of an infra-nilmanifold 
endomorphism as given by M.W.~Hirsch is the best possible. In fact, we will show that the only 
maps of an infra-nilmanifold, that lift to an automorphism of the corresponding nilpotent Lie group are 
exactly the infra-nilmanifold endomorphisms defined in Definition~\ref{defHirsch}. When reading the work 
of J.~Franks, it is clear that he also only considers those maps on an infra-nilmanifold $\Gamma\backslash N$ which 
lift to an automorphism of the Lie group $N$ (E.g.~see the first few lines of the proof of Theorem~2.2 of \cite{fran70-1}). 
In fact, when talking about infra-nilmanifold endomorphisms most authors, including J.~Franks, M.~Schub and M.~Hirsch (but e.g.\
also in \cite{bl93-1,hass02-1,hira88-1,sumi96-1} and in many others papers) are talking about maps which lift to an automorphism 
of the Lie group $N$. 

\begin{theorem}\label{hirschok}
Let $N$ be a connected and simply connected nilpotent Lie group, $F\subseteq \Aut(N)$ a finite group and 
$\Gamma$ a torsion free discrete and uniform subgroup of $N\semi F$ and assume that the holonomy 
group of $\Gamma$ is $F$. If $\varphi:N\rightarrow N$ 
is an automorphism for which the map 
\[ \bar\varphi:\; \Gamma\backslash N \rightarrow \Gamma \backslash N:\; \Gamma\cdot  n \mapsto \Gamma\cdot  \varphi(n) \]
is well defined (meaning that $\Gamma\cdot  \varphi(n) = \Gamma\cdot  \varphi (\gamma\cdot n)$ for all $\gamma\in \Gamma$), 
then 
\[\Phi : N \semi F \rightarrow N \semi F :\; x \mapsto \varphi x \varphi^{-1} \mbox{ (conjugation in $\Aff(N)$)}\]
is an automorphism of $N\semi F$, with $\Phi(F)=F$ and $\Phi(\Gamma)\subseteq \Gamma$. Hence,
$\bar\varphi$ is a infra-nilmanifold endomorphism (as in Definition~\ref{defHirsch}). 
\end{theorem}

\noindent\underline{Proof:}
The fact that $\bar\varphi$ is well defined, means that $\varphi$ is a lift of $\bar\varphi$ to the universal cover 
$N$ of $\Gamma\backslash N$. Now, $\forall \gamma \in \Gamma$, also the composition $\varphi \gamma$ is a lift of $\bar
\varphi$, since $\Gamma$ is the group of covering transformations of the covering $N\rightarrow \Gamma\backslash N$. 
It follows that there exists a $\gamma'$ such that $\varphi \gamma = \gamma' \varphi$. Now, since $\varphi$ is an automorphism
of $N$, we can write this as $\varphi \gamma \varphi^{-1} = \gamma'$ for some $\gamma'\in \Gamma$ so 
\[ \varphi \Gamma \varphi^{-1} \subseteq \Gamma.\]
Now, consider the inner automorphism $\Psi$ of $\Aff(N)$ induced by $\varphi$:
\[ \Psi:\;\Aff(N)\rightarrow \Aff(N):\; x \mapsto \varphi x \varphi^{-1}.\]
For all $n\in N$, we have that $\Psi(n)= \varphi(n)$, so $\Psi(N)=\varphi(N)=N$. We showed above that that 
$\Psi(\Gamma)\subseteq \Gamma$. It follows that $\Psi(N\semi F)= \Psi(N\Gamma)\subseteq N\Gamma$. Hence,
$\Psi$ induces an injective endomorphism of $N\semi F$. As $F$ is mapped into itself by $\Psi$ (because
$\Aut(N)$ is mapped into itself by $\Psi$) and $F$ is finite, we must have that $\Psi(F)=F$. Together with the
fact that $\Psi(N)=N$, this implies that $\Psi(N\semi F)=N\semi F$ and hence
$\Psi$ restricts to an automorphism $\Phi$ of $N\semi F$, satisfying the conditions mentioned in the statement of the 
theorem. \qed

\bigskip

\begin{remark} When checking literature, it seems that most authors that are talking about infra-nilmanifold
endomorphisms, seem to assume that such a map lifts to an automorphism of the covering Lie group $N$. Hence, this
implies that they are actually using the definition of M.W.Hirsch (which is probably also 
the definition that J.~Franks meant to give). So from now onwards, when we use the term infra-nilmanifold 
endomorphism, we are referring to the only correct Definition~\ref{defHirsch}. 
\end{remark}

\bigskip

We are now ready to define the generalization of the concept of an infra-nilmanifold endomorphism we announced in the
introduction.

\begin{definition}
Let $N$ be a connected and simply connected nilpotent Lie group, $F\subseteq \Aut(N)$ a finite 
group, $\Gamma$ a torsion free discrete and uniform subgroup of $N\semi F$. 
Let $\alpha\in \Aff(N)$ be an element such that $\alpha \Gamma \alpha^{-1} \subseteq \Gamma$, 
then $\alpha$ induces a map 
\[ \bar{\alpha} : \Gamma\backslash N \rightarrow \Gamma \backslash N:\; \Gamma\cdot n \mapsto \Gamma\cdot \alpha(n).\]
We call $\bar{\alpha}$ an affine endomorphism of the infra-nilmanifold $\Gamma \backslash N$ induced by $\alpha$. 
When 
$\alpha \Gamma \alpha^{-1} =\Gamma$, the map $\bar\alpha$ is a diffeomorphism, and we call 
$\bar\alpha$ an affine automorphism.
\end{definition}

As it is so crucial for what follows, we briefly recall from the theory of covering transformations
how the group $\Gamma$ can be seen as the fundamental group of $\Gamma$ and what the effect of an affine 
endomorphism is on the fundamental group. Details of what follows can be found in any text book dealing with 
this topic, e.g.~\cite[Chapter 2]{span66-1} and \cite[Chapter 5]{mass67-1}.

Choose any basepoint $n_0\in \Gamma\backslash N$ and choose a point $\tilde{n}_0\in N$ whose orbit corresponds to the 
point $n_0$. Now, any loop $f:I\rightarrow \Gamma\backslash N$ at $n_0$ ($I$ is the unit interval $[0,1]$) has a unique lift 
to a path $\tilde{f}:I \rightarrow N$ starting at $\tilde{n}_0$ (i.e\ $\tilde{f}(0)=\tilde{n}_0$).
The endpoint $\tilde{n}_1=\tilde{f}(1)$ of $\tilde f$ lies in the same orbit as $\tilde{n}_0$ (because they both project 
onto $n_0$) and hence, there exists a $\gamma_f \in \Gamma$ with $\gamma_f\cdot \tilde{n}_0 = \tilde {n_1}$. In this way, we associate 
to any loop  $f$ at $n_0$ an element $\gamma_f\in \Gamma$. 
It is a general fact that this correspondence does not depend on the path homotopy class of $f$ and 
defines an isomorphism $\Phi:\Pi_1(\Gamma\backslash N, n_0) \rightarrow \Gamma$. Note that this 
isomorphism depends on the choice of the point $\tilde{n}_0$ and that a different choice, say $\tilde{n}_1$, 
changes the isomorphism by an inner automorphism of $\Gamma$.

\medskip

Now, let $\bar\alpha$ be an affine endomorphism induced by an affine map $\alpha \in \Aff(N)$ (with 
$\alpha \Gamma \alpha^{-1} \subseteq \Gamma$).
Choose a basepoint $n_0\in \Gamma\backslash N$ and a point $\tilde{n}_0\in N$ projecting onto $n_0$.
Then $\tilde{n}_1= \alpha(\tilde{n}_0)\in N$ is a point projecting onto $n_1=\alpha(n_0)$. Now, let us use 
$\tilde{n}_0$ resp.\ $\tilde{n}_1$ to identify $\Pi_1(\Gamma\backslash N, n_0)$ resp.\  $\Pi_1(\Gamma\backslash N, n_1)$ 
with $\Gamma$. Let $\bar\alpha_\sharp:\Pi_1(\Gamma\backslash N, n_0)\rightarrow \Pi_1(\Gamma\backslash N, n_1)$ denote 
the morphism induced by $\bar\alpha$. We claim that $\bar\alpha$ is exactly conjugation with $\alpha$. Indeed,
consider again a loop $f$ based at $n_0$ and let $\tilde{f}$ be the lift of $f$ to $N$ starting at $\tilde{n}_0$.
Let $\gamma\in \Gamma$ be the element such that $\gamma\cdot \tilde{n}_0$ is the endpoint of $\tilde{f}$ (so the path class 
$[f]\in \Pi_1(\Gamma\backslash N, n_0)$ corresponds to $\gamma\in \Gamma$). It is obvious that 
$\alpha \circ \tilde{f}$ is the unique lift, beginning in $\tilde{n}_1$,  of the loop $\bar\alpha \circ f$.
The endpoint of $\alpha \circ \tilde{f}$ is 
\[ \alpha (f(1))= \alpha(\gamma\cdot \tilde{n}_0) = \alpha (\gamma (\alpha^{-1} ( \alpha(  \tilde{n}_0)))) = (\alpha \gamma \alpha^{-1}) \cdot \tilde{n}_1.\]
This shows that the element of $\Gamma$ corresponding to $\bar\alpha_\sharp [f]= [\bar\alpha \circ f]$ is exactly $\alpha \gamma \alpha^{-1}$.

\medskip

Note that in the discussion above, we have chosen $\tilde{n}_1$ based on our knowledge of $\alpha$. In practice, this is often not possible or even not desirable.
E.g.\ in this paper we often choose a fixed point $n_0$ of a selfmap $\bar\alpha$ on an infra-nilmanifold as a base point. To study then 
the induced morphism $\bar\alpha_\sharp: \Pi_1(\Gamma\backslash N, n_0) \rightarrow  \Pi_1(\Gamma\backslash N, n_0)$ we will of course use 
two times the same $\tilde{n}_0$ when identifying $\Pi_1(\Gamma\backslash N, n_0)$ with $\Gamma$. This implies that $\bar\alpha_\sharp$ will 
only be the same as conjugation with $\alpha$ in $\Aff(N)$ up to an inner conjugation by an element of $\Gamma$.

\medskip

It follows that we have the following 

\begin{proposition}\label{laternodig}
Let  $\bar\alpha$ be an affine endomorphism of an infra-nilmanifold $\Gamma\backslash N$ and let 
$\psi: \Gamma\rightarrow \Gamma: \gamma \mapsto \alpha \gamma \alpha^{-1}$ be the corresponding 
monomorphism of $\Gamma$. Then, the map $\bar\alpha_\sharp:\Gamma=\Pi_1(\Gamma\backslash N,x) \rightarrow
\Gamma=\Pi_1(\Gamma\backslash N,\bar{\alpha}(x))$ is, up to composition with an inner automorphism of $\Gamma$, precisely $\psi$. 
\end{proposition}

\begin{remark}\label{model-prop}
At this point, it is worthwhile to indicate that \cite[Proposition 3.5]{fran70-1}, which is crucially used at other places in the work of Franks 
(e.g.~in the basis theorem \cite[Theorem 8.2]{fran70-1} on which Gromov's result is based), is not correct. This Proposition claims that for 
any covering $f:K\rightarrow K$, where $\Pi_1(K)$ is a finitely generated, torsion free and virtually nilpotent group, there exists
an infra-nilmanifold $M$ and an infra-nilmanifold endomorphism $g:M\rightarrow M$ which is $\Pi_1$--conjugate to $f$. 
This is not true (see the example below and the examples in the following sections) and one really also needs to consider affine endomorphisms
of the infra-nilmanifolds as well. On the other hand, when $\Pi_1(K)$ is nilpotent (or abelian) the proposition is correct. 

The problem in the alleged proof is situated at the very end of it on page 78. First of all, the wrong result of Auslander is used (but 
this can be solved by using Corollary~\ref{correctAusl}). However, as indicated by the example above, the automorphism $\bar{g}$ (where 
I now use the notations of \cite[page 78]{fran70-1}) does not necessarily induce a map on the infra-nilmanifold $M$ (and even if it does,
the induced map on the fundamental group is not necessarily the map $g_\ast$). 
\end{remark}

We finish this section by giving a counter-example to Fanks' ``Existence of a Model''--proposition (\cite[Proposition 3.5]{fran70-1}).
Consider the Klein Bottle $K$ and choose a base point $x_0\in K$. Then, the fundamental group $\Pi_1(K,x_0)\cong 
\Gamma=\langle a,b \;|\; ba = a^{-1} b \rangle$. It is easy to find a homeomorphism $f:K\rightarrow K$, with $f(x_0)=x_0$ and such that 
$f_\sharp: \Pi_1(K,x_0) \rightarrow \Pi_1(K,x_0)$ satisfies $f_\sharp (a) = a$ and $f_\sharp(b)=ab$. Now, consider any embedding of $\Gamma$ 
into $\Isom(\R^2)$ as a discrete subgroup, then the translation subgroup of $\Gamma$ will be $\Gamma\cap \R^2=\langle a, b^2\rangle\cong \Z^2$.
Now, assume that $\bar\varphi: \Gamma\backslash \R^2 \rightarrow \Gamma\backslash \R^2$ is an infra-nilmanifold endomorphism (induced 
by the automorphism $\varphi:\R^2 \rightarrow \R^2$), which is $\Pi_1$--conjugate to $f$. This means that there is a commutative diagram
\[ \xymatrix{
\Pi_1(K,x_0)\ar[r]^-{\Phi}\ar[d]_{f_\sharp} & \Pi_1(\Gamma\backslash \R^2, \Gamma\cdot 0)\cong \Gamma\ar[d]^{\bar\varphi_\sharp} \\
\Pi_1(K,x_0)\ar[r]_-{\Phi}                  & \Pi_1(\Gamma\backslash \R^2, \Gamma\cdot 0)\cong \Gamma       \\
}\] 
for some isomorphism $\Phi$. As $f_\sharp(a)= a$ and $f_\sharp(b^2)= b^2$, it follows that $\bar\varphi_\sharp$ has to be the identity 
on the translation subgroup $\langle a, b^2\rangle$ of $\Gamma$. But as the restriction of $\bar\varphi_\sharp$ to the translation subgroup 
is exactly the same as the restriction of $\varphi$ to this translation subgroup, it follows that $\varphi$ is the identity on this translation subgroup 
and hence $\varphi$ is just the identity automorphism of $\R^2$. But this means that $\bar\varphi$ is the identity map also, hence $\bar\varphi_\sharp$ 
is the identity automorphism, which contradicts the commutativity of the diagram above. 

\section{An expanding map not topologically conjugate to an in\-fra-nilmanifold endomorphism}
Already on the smallest example of an infra-nilmanifold which is not a nilmanifold (or a torus) we can
construct an expanding map which is not topologically conjugate to an infra-nilmanifold endomorphism of that infra-nilmanifold.
Our example will be an affine endomorphism of the Klein Bottle.
This example shows that there are problems with the proof of the geometric corollary 
on page 55 of \cite{grom81-1}, which we will explain below. 
Of course, this does not cast any doubt on the (very nice) main result of \cite{grom81-1}
stating that finitely generated groups of polynomial growth are virtually nilpotent!

\bigskip

So, consider the Klein Bottle which is constructed by taking the group $\Gamma\subseteq \R^2\semi \Z_2$,
where 
\begin{equation}\label{Z2}
\Z_2=\left\{ \left(\begin{array}{cc} 1 & 0 \\ 0 & 1 \end{array} \right),\;
               \left(\begin{array}{cc} -1 & 0 \\ 0 & 1 \end{array} \right)\right\}\subseteq \GL_2(\R).
\end{equation}
The torsion free discrete and uniform subgroup $\Gamma$ of $\R^2\semi \Z_2$ we use to construct the Klein 
Bottle is generated by the following 2 elements:
\[ a=\left( \left( \begin{array}{c}1 \\  0  \end{array}\right),\;
 \left(\begin{array}{cc} 1 & 0 \\ 0 & 1 \end{array} \right) \right)\mbox{ and }
b=\left( \left( \begin{array}{c} 0 \\ \frac12   \end{array}\right),\;
 \left(\begin{array}{cc} -1 & 0 \\ 0 & 1 \end{array} \right) \right)\]
Note that $a$ and $b^2$ generate the group of translations $\Z ^2$. 
Let $\alpha$ be the affine map
\[\alpha:\R^2\rightarrow \R^2: \left( \begin{array}{c} x \\ y \end{array}\right) \mapsto 
\left( \begin{array}{c} 3x+\frac12 \\ 3 y \end{array}\right).\]
So 
\[\alpha = \left(  \left( \begin{array}{c}\frac12 \\  0  \end{array}\right),\;
 \left(\begin{array}{cc} 3 & 0 \\ 0 & 3 \end{array} \right) \right)\in \Aff(\R^2).\]
One easily checks that 
\[ \alpha a \alpha^{-1} = a^3 \mbox{ and } \alpha b \alpha^{-1} = a b^3\]
showing that $\alpha \Gamma \alpha^{-1} \subseteq \Gamma.$ Hence $\alpha$ induces an 
affine endomorphism $\bar\alpha:\Gamma\backslash \R^2 \rightarrow \Gamma \backslash \R^2$ of 
the Klein bottle $K=\Gamma \backslash \R^2$. Moreover, as the linear part of $\alpha$ has only eigenvalues 
of modulus $>1$, the map $\bar \alpha$ is an expanding map of the Klein bottle.

I claim that this map is not topologically conjugate to an expanding infra-nilmanifold endomorphism of this 
Klein Bottle.

To see this, 
suppose that $\varphi:\R^2 \rightarrow \R^2$ is a linear isomorphism  inducing an endomorphism 
$\bar\varphi:\Gamma\backslash \R^2 \rightarrow \Gamma \backslash \R^2$ of the Klein bottle. By 
Theorem~\ref{hirschok}, we know that 
$\varphi \Gamma \varphi^{-1} \subseteq \Gamma$. From this, it also follows that 
$\varphi \Z_2 \varphi^{-1}= \Z_2$, where $\Z_2\subseteq \GL(2,\R)$ is  as in (\ref{Z2}).
Hence,
\[ \varphi \left(\begin{array}{cc} -1 & 0 \\
0 & 1 \end{array}\right) = \left(\begin{array}{cc} -1 & 0 \\
0 & 1 \end{array}\right) \varphi, \] 
from which it follows that 
\[ \varphi=\left(\begin{array}{cc} k & 0 \\
0 & l \end{array}\right)\]
for some $k,l\in \R$.
Now, requiring that $\varphi a \varphi^{-1} \in \Gamma$ and $\varphi b \varphi^{-1} \in \Gamma$ 
leads to the condition that $k\in \Z$ and $l=2 m+1$, for $m\in \Z$ (so $l$ is odd).

\medskip

As recalled in some detail in the discussion before Proposition~\ref{laternodig}, there is 
an isomorphism $\Pi_1(\Gamma\backslash \R^2, \bar{0})\cong\Gamma$. Here, we use $\bar{0}$ to denote the image of the zero vector 
in the Klein Bottle $\Gamma\backslash \R^2$ and we use the zero vector as the point $\tilde{n}_0$ (see discussion before Prop.~\ref{laternodig})
to establish the isomorphism between $ \Pi_1(\Gamma\backslash \R^2, \bar{0})$ and $\Gamma$. From proposition~\ref{laternodig}, we know that 
the map induced by $\bar\varphi$ is the same as conjugation with $\varphi$ inside $\Aff(\R^2)$.

Now, suppose that $\bar\alpha$ is topologically conjugate to $\bar\varphi$, then 
there must exist a homeomorphism $h:\Gamma\backslash \R^2 \rightarrow \Gamma\backslash \R^2 $, such that $h \circ \bar\alpha = \varphi \circ h$.
Now, choose $h^{-1}(\bar{0})$ as another basepoint of $\Gamma\backslash \R^2$. It is obvious that $h^{-1}(\bar{0})$ is a fixed point of 
$\bar\alpha$. We know that we can also fix an isomorphism of $\Pi_1(\Gamma\backslash \R^2, h^{-1}(\bar{0}))$ with $\Gamma$ and that 
under this identification the map $\bar{\alpha}_\sharp:\Gamma\rightarrow \Gamma$ is , up to an inner automorphism, exactly the same as 
conjugation with $\alpha \in \Aff(\R^2)$. 

\medskip

\noindent Using the above, we find a commutative diagram of groups and morphisms
\[ \xymatrix{ \Gamma \ar[r]^{h_\sharp} \ar[d]_{\bar{\alpha}_\sharp}  & \Gamma \ar[d]^{\bar{\varphi}_\sharp} \\
\Gamma  \ar[r]_{h_\sharp}  & \Gamma
}
\]
This diagram leads to an induced diagram of morphisms on the abelianization of $\Gamma$:
\[ \xymatrix{ \Gamma/[\Gamma,\Gamma] \ar[r]^{h_\ast} \ar[d]_{\bar{\alpha}_\ast}  & \Gamma/[\Gamma,\Gamma] \ar[d]^{\bar{\varphi}_\ast} \\
\Gamma/[\Gamma,\Gamma]  \ar[r]_{h_\ast}  & \Gamma/[\Gamma,\Gamma]
}
\]
We have that 
$\Gamma/[\Gamma,\Gamma] =\Z_2 \oplus \Z$, where $\Z_2$ (resp.~$\Z$) is generated by 
the natural projection $\bar a$ of $a$ (resp.~$\bar b$ of $b$). 

\medskip

\noindent As $\bar\alpha_\sharp$ was, up to an inner automorphism of $\Gamma$, the same as conjugation with $\alpha$ inside $\Aff(\R^2)$,
we know exactly what $\bar\alpha_\ast$ is, and we also already obtained some information on $\varphi_\ast$:
\[ \bar\alpha_\ast (\bar a) = \bar a^3=\bar a,\; \bar\alpha_\ast (\bar b)= \bar a\bar b^3, \;
\bar\varphi_\ast (\bar a) = \bar a^l,\; \bar\varphi_\ast (\bar b)= \bar b^{2 m+1} \mbox{ with }l,m\in \Z. \] 
As $h$ is a homeomorphism of the Klein bottle, we know that $h_\ast$ is an isomorphism of 
$\Gamma/[\Gamma,\Gamma]$. It follows that $h_\ast(\bar a) = \bar a$ while 
for $h_\ast (\bar b)$ we have one of the following four possibilities:
\[ h_\ast (\bar b) = \bar b,\; 
h_\ast (\bar b) = \bar b^{-1},\;
h_\ast (\bar b) = \bar a\bar b \mbox{ or } 
h_\ast (\bar b) = \bar a\bar b^{-1}\]
It is now easy to see that for none of these four possibilities, we can have that 
\[ h_\ast \circ \bar\alpha_\ast = \bar\varphi_\ast \circ h_\ast,\]
contradicting the fact that $h \circ \bar\alpha = \varphi \circ h$ and hence showing that 
$\bar\varphi$ is not topologically conjugate to $\bar\alpha$. \qed

\medskip

This example indicates a real problem in the proof of the geometric corollary on page 55 of \cite{grom81-1}. 
In fact, this geometric corollary follows from Gromov's main result by applying \cite[Theorem 8.3]{fran70-1} 
and \cite[Theorem 5]{shub70-1} (or the equivalent \cite[Theorem 8.2]{fran70-1}). Now looking at the proof 
of \cite[Theorem 5]{shub70-1} (or \cite[Theorem 8.2]{fran70-1}) ones sees that actually the incorrect 
``Existence of a Model''--Proposition of Franks is used (see remark~\ref{model-prop}).  In fact, both Shub and 
Franks are claiming that an expanding map on an infra-nilmanifold is topologically conjugate to an expanding 
infra-nilmanifold endomorphism of the same infra-nilmanfiold, which is actually wrong by the example above.

\medskip\

However, in the sequel of this section, we will show that any expanding map of a given infra-nilmanifold 
is topologically conjugate to an expanding affine endomorphism of the same infra-nilmanifold, from which it will 
follow that any expanding map of a compact manifold $M$ will be topologically conjugate to an expanding 
affine infra-nilmanifold endomorphism of any infra-nilmanifold with the same fundamental group as $M$.

\medskip

In order to prove this result, we need some more results concerning affine maps of 
infra-nilmanifolds. Let $N$ be a connected and simply connected nilpotent Lie group,
$\delta\in \Aut(N)$ and $d\in N$. Then $D=(d,\delta)$ is an affine map of 
$N$. As $\delta\in \Aut(N)$, we know that its differential $\delta_\ast\in \Aut(\lien)$,
where $\lien$ is the Lie algebra of $N$. When we talk about the eigenvalues of $D$ (or the 
eigenvalues of $\delta$) we will mean the eigenvalues of $\delta_\ast$. 

\begin{lemma}\label{eigenwaarde1}
Let $N$ be a connected and simply connected nilpotent Lie group and $D\in \Aff(N)$.
If 1 is not an eigenvalue of $D$, then there is a unique fixed point $n_0\in N$ for 
the affine map $D$. 
\end{lemma}
\underline{Proof:} This is a special case of \cite[Lemma 2]{bdd04-1}. \qed

\medskip

Now, consider a finitely generated and torsion free nilpotent group $\Lambda$ and an 
injective endomorphism $\varphi\in \Aut(\Lambda)$. Up to isomorphism there is a unique connected 
and simply connected nilpotent Lie group $N$,  containing $\Lambda$ as a uniform discrete 
subgroup. This $N$ is called the Mal'cev completion of $\Lambda$. The endomorphism 
$\varphi$ extends uniquely to a continuous automorphism $\tilde\varphi\in \Aut(N)$ and we 
can talk about the eigenvalues of $\varphi$, by which we will mean the eigenvalues of 
$\tilde\varphi$ (which in their turn are the eigenvalues of the differential 
$\tilde{\varphi}_\ast\in \Aut(\lien)$ of $\tilde{\varphi})$.

\medskip

More generally, we can consider as before a torsion free uniform discrete 
subgroup $\Gamma\subseteq N\semi F$, where $N$ is a connected and simply connected 
nilpotent Lie group and $F$ is a finite subgroup of $\Aut(N)$. We assume that $F$ is the 
holonomy group of $\Gamma$. We know that $\Lambda=\Gamma\cap N$
is a uniform discrete subgroup of $N$ and so $N$ is the Mal'cev completion of $\Lambda$.
Let $\varphi:\Gamma \rightarrow \Gamma $ be an injective endomorphism of $\Gamma $. It follows from 
Corollary~\ref{correctAusl} that $\varphi$ extends uniquely to an automorphism of $N\semi F$ and restricts 
to an injective endomorphism of $\Lambda$. We define the eigenvalues of $\varphi$ to be the 
eigenvalues of the restriction of $\varphi$ to $\Lambda$. \\
On the other hand, we know that $\varphi$ can also be realized as conjugation by some element 
$D=(d,\delta)$ in $\Aff(N)$. It turns out that the eigenvalues of $\varphi$ are exactly the same as the eigenvalues 
of $D$.

\begin{lemma}\label{zelfdeeig}
 Let $N$ be a connected and simply connected nilpotent Lie group and 
Let $F$ be a finite subgroup of $\Aut(N)$. Assume that  
$\Gamma$ is  a uniform discrete subgroup of $N\semi F$ with holonomy group $F$,  $\varphi$ is 
an injective endomorphism of $\Gamma$ and that $D=(d,\delta)\in \Aff(N)$ realizes this endomorphism via 
conjugation in $\Aff(N)$:
\[ \forall \gamma\in \Gamma:\; \varphi(\gamma) = (d,\delta) \gamma (d,\delta)^{-1}.\]
Then, the set of eigenvalues of $\varphi$ is exactly the same as the set of eigenvalues 
of $D$.
\end{lemma}

\underline{Proof:} To  compute the eigenvalues of $\varphi$, we have to find the eigenvalues 
of the induced automorphism $\tilde\varphi$ of $N$ (obtained by first extending $\varphi$ to $N\semi F$ and then taking the 
restriction to $N$).
But this automorphism is also obtained by conjugation with $D$:
\[ \forall n\in N:\; \tilde\varphi (n)= DnD^{-1} = (d,\delta) n (d,\delta)^{-1} = d \delta(n) d^{-1}= (\mu(d)\circ \delta)(n).\]
where $\mu(d)$ denotes conjugation with $d\in N$.
It follows that the eigenvalues of $\varphi$ are precisely the same as the eigenvalues of $\mu(d)\circ \delta$. It is a standard
argument to show that an inner automorphism of a nilpotent Lie group has no influence on the eigenvalues: indeed,
to find the eigenvalues of a given automorphism $\psi\in \Aut(N)$, we can consider the filtration of $N$ by the terms of its 
lower central series (which goes to 1 as $N$ is nilpotent)
\[ N=\gamma_1(N)\supseteq \gamma_2(N)\supseteq \cdots \supseteq \gamma_i(N)\supseteq \gamma_{i+1}(N)=[\gamma_i(N),N]\supseteq \cdots \supseteq \gamma_c(N)=1.\]
Each term in this filtration is invariant under $\psi$ and analogously the corresponding terms of the lower central series of 
the Lie algebra $\lien$ of $N$:
\[ \lien=\gamma_1(\lien)\supseteq \gamma_2(\lien) \supseteq \cdots \supseteq \gamma_i(\lien)\supseteq 
\gamma_{i+1}(\lien)=[\gamma_i(\lien),\lien]\supseteq \cdots \supseteq \gamma_c(\lien)=1\]
are then invariant under the differential $\psi_\ast$ of $\psi$. It follows that to find the eigenvalues of $\psi$, we have to find 
the eigenvalues of the induced automorphism on each quotient $\gamma_i(\lien)/\gamma_{i+1}(\lien)$. However, an inner automorphism 
of $N$ induces the identity on each quotient $\gamma_i(N)/\gamma_{i+1}(N)$ and so its differential induces the identity on 
$\gamma_i(\lien)/\gamma_{i+1}(\lien)$. It follows that $\delta$ and $\mu(d)\circ \delta$ induce the same linear map on each 
quotient $\gamma_i(\lien)/\gamma_{i+1}(\lien)$ and hence, they have the same eigenvalues. \qed

\medskip

In what follows it will be crucial to know when an affine map does not have 1 as an eigenvalue (so that we 
will be able to apply Lemma~\ref{eigenwaarde1}). The following lemma can serve as a criterion for this.

\begin{lemma} \label{crit1}
Let $N$ be a connected and simply connected nilpotent Lie group and 
Let $F$ be a finite subgroup of $\Aut(N)$. Assume that  
$\Gamma$ is  a uniform discrete subgroup of $N\semi F$ with holonomy group $F$ and $\varphi$ is 
an injective endomorphism of $\Gamma$. 
If $\varphi$ has 1 as an eigenvalue, then there exists a non-trivial subgroup $\Delta$ 
of $\Gamma$ such that for all $\gamma\in \Delta:\; \varphi(\gamma)=\gamma$.
\end{lemma}

\underline{Proof:}
Let $\Lambda=\Gamma\cap N$. As already argued above, $\varphi$ restricts to an injective endomorphism of 
$\Lambda$ and this restriction extends uniquely to an automorphism of $N$. We will use the same symbol $\varphi$ 
to denote all these endomorphisms. Recall that for nilpotent Lie groups, the 
exponential map $\exp: \lien \rightarrow N $ is a diffeomorphism ($\lien$ is the Lie algebra of $N$) and we denote 
its inverse by $\log$. Consider now $\lien_\Q= \Q \log(\Lambda)$  (the rational span of $\log(\Lambda)$) and $N_\Q=\exp (\lien_\Q)$. 
The vector space $\lien_\Q$ is a rational Lie algebra and the differential $\varphi_\ast$ of $\varphi$ restricts to an automorphism 
of $\lien_\Q$. For more details about this and following facts on these rational Lie algebras, we refer to \cite[Chapter 6]{sega83-1}. 
As $\psi_\ast$ has 1 as an eigenvalue, there exists a nonzero vector $X\in \lien_\Q$ with $\varphi_\ast(X)=X$. This implies 
that $1\neq x=\exp(X)\in N_\Q$ is an element with $\varphi(x)=x$. Now, $N_\Q$ is the radicable hull (\cite[Page 107]{sega83-1}) of $\Lambda$, and so there 
exists a positive integer $k>0$ such that $1\neq x^k \in \Lambda$. It follows that $x^k$ is a nontrivial element of $\Lambda$ 
with $\varphi(x^k)= x^k$. The proof now finishes by taking $\Delta$ to be the group generated by $x^k$. \qed

\medskip

We are now ready to prove the main result of this paper in which we will adopt J.~Franks' original approach 
for infra-nilmanifold endomorphisms \cite[Section 8]{fran70-1} to  the more general case of affine endomorphisms.

\begin{theorem}\label{main1}
Let $f:M\rightarrow M$ be an expanding map of a compact manifold $M$. Then, there exists an infra-nilmanifold 
$\Gamma\backslash N$ whose fundamental group $\Gamma$ is isomorphic to $\Pi_1(M)$. And for any such $\Gamma\backslash N$, there exists 
an expanding affine endomorphism of that infra-nilmanifold which is topologically conjugate to $f$. 
\end{theorem}

\underline{Proof:}
By \cite[Theorem 1]{shub69-1} we can choose a fixed point $m_0\in M$  of $f$.  From \cite[Theorem 8.3]{fran70-1} we know that 
$\Pi_1(M,m_0)$ has polynomial growth and so by the main result of \cite{grom81-1} it follows that $\Pi_1(M,m_0)$ has a nilpotent subgroup
of finite index.  Moreover, by \cite[Proposition 3]{shub69-1}, we know that 
$\Pi_1(M,m_0)$ is torsion free, $M$ is a $K(\Pi_1(M,m_0),1)$--space  and the induced map $f_\sharp:\Pi_1(M,m_0) \rightarrow \Pi_1(M,m_0)$ is an injective endomorphism.
Every finitely generated virtually nilpotent group can be realized as a uniform and discrete subgroup of a semi-direct product 
$N\semi F$, where $N$ is a connected and simply connected nilpotent Lie group and $F$ is a finite subgroup of $\Aut(N)$ (e.g~ \cite[Theorem 3.1.3]{deki96-1}).
Fix such an embedding $i:\Pi_1(M,m_0) \rightarrow N\semi F$ realizing $\Pi_1(M,m_0)$ as such a uniform discrete  subgroup and denote $\Gamma= i (\Pi_1(M,m_0))$. 
Without loss of generality we assume that $F$ is the holonomy group of $\Gamma$. 
So there is an isomorphism $A:\Gamma \rightarrow \Pi_1(M,m_0)$ (where $A$ is in fact the inverse of $i$), already showing the existence of the 
infra-nilmanifold $\Gamma\backslash N$.  We continue our proof with a fixed choice of such an infra-nilmanifold. Let $B= A^{-1} \circ f_\sharp \circ A$, then $B$ 
is an injective endomorphism of $\Gamma$ and so there exists an affine map $\alpha=(d,\delta)\in \Aff(N)$ with 
$B(\gamma) = \alpha \gamma \alpha^{-1}$, for all $\gamma\in \Gamma$.
By \cite[Corollary 1]{shub69-1} we know that the identity element is the unique fixed element of $f_\sharp$ and so the identity element is also 
the only fixed point for $B$.  By Lemma~\ref{crit1} it follows that $\alpha$ does not have 1 as one of its eigenvalues and so, by Lemma~\ref{eigenwaarde1} there 
exists a unique fixed point $\tilde{n}_0\in N$ for $\alpha$. Let $n_0$ be the corresponding point in the infra-nilmanifold $\Gamma\backslash N$ 
and use the point $\tilde{n}_0$ to identify the fundamental group $\Pi_1(\Gamma\backslash N,n_0)$ with $\Gamma$. By the discussion before Proposition~\ref{laternodig},
we know that $\alpha$ induces an affine endomorphism $\bar\alpha$ of $\Gamma\backslash N$, with $n_0$ as a fixed point, and that the induced endomorphism $\bar\alpha_\sharp$ 
of $\Pi_1(\Gamma\backslash N,n_0)=\Gamma$ is exactly $B$. We therefore have a commutative diagram
\[ 
\xymatrix{\Pi_1(\Gamma\backslash N,n_0)=\Gamma \ar[r]^{\;\;\;A} \ar[d]_{\bar\alpha_\sharp} &  \Pi_1(M,m_0)\ar[d]^{f_\sharp} \\
          \Pi_1(\Gamma\backslash N,n_0)=\Gamma \ar[r]_{\;\;\;A}  &  \Pi_1(M,m_0)}
\]
By \cite[Theorem 4]{shub69-1} there exists a unique continuous map $h:(\Gamma\backslash N,n_0) \rightarrow (M,m_0) $ with $f \circ h = h \circ \bar\alpha$ 
and for which $h_\sharp : \Pi_1(\Gamma\backslash N,n_0) \rightarrow \Pi_1(M,m_0)$ is exactly $A$.  (As usual, by a map $g:(X,x)\rightarrow (Y,y)$ we mean a map from 
the space $X$ to the space $Y$, with $g(x)=y$ where $x$ and $y$ are given points of $X$ and $Y$ respectively). As $A$ is an isomorphism, $h$ is a homotopy equivalence, since
we are working with $K(\Pi,1)$--spaces.

Let $\tilde{M}$ denote the universal covering space of $M$, with covering projection $p_M:\tilde{M}\rightarrow M$ and let $\tilde{m}_0\in \tilde{M}$ be 
a point with $p_M(\tilde{m}_0) = m_0$.  Now, consider the unique lift $\tilde{h}:(N,\tilde{n}_0) \rightarrow (\tilde{M},\tilde{m}_0)$ of $h$ and the 
unique lift $\tilde{f}:(\tilde{M},\tilde{m}_0) \rightarrow (\tilde{M},\tilde{m}_0)$ of $f$, then 
\[ \tilde{f} \circ \tilde{h} = \tilde{h} \circ \alpha.\]
Let $L_{\tilde{n}_0}:N\rightarrow N: x \mapsto \tilde{n}_0\, x$ denote left translation by $\tilde{n}_0$ in $N$. As $\tilde{n}_0$ is a fixed point of $\alpha=(d,\delta)$, 
we have that $L_{\tilde{n}_0} \circ \delta= \alpha \circ L_{\tilde{n}_0}$. Summarizing the above, we obtain the following commutative diagram of maps and spaces in 
which $\exp$ and $L_{\tilde{n}_0}$ are diffeomorphisms.
\[
\xymatrix{
\lien\ar[d]_{\delta_\ast}\ar[r]^\exp  &    N \ar[d]_\delta \ar[r]^{L_{\tilde{n}_0}} & N\ar[d]_\alpha\ar[r]^{\tilde{h}}   & \tilde{M}\ar[d]_{\tilde{f}}\\
\lien \ar[r]_\exp  &    N \ar[r]_{L_{\tilde{n}_0}}  & N\ar[r]_{\tilde{h}}   & \tilde{M}\\
}
\]
Let $k= \tilde{h} \circ L_{\tilde{n}_0} \circ \exp$. By \cite[Lemma 3.4]{fran70-1}, the map $\tilde{h}$ and hence also $k$ is a proper map.
We can now continue as in Franks' paper to show that $\delta_\ast$ only has eigenvalues of modulus $>1$.
From $\tilde{f} \circ k = k \circ \delta_\ast$ it immediately follows that $\tilde{f}^n \circ k = k \circ \delta_\ast^n$. Now, assume that 
$\delta_\ast$ has an eigenvalue of modulus $\leq 1$. It then follows that there exists a non-zero element $x\in \lien$ with 
$\| \delta_\ast^n(x)\| \leq \| x \|$ (where $\| \;\;\|$ denotes the usual norm on $\lien$). (Note that the argument given in \cite{fran70-1} is not 
completely correct, because he considers an eigenvector of the corresponding eigenvalue of modulus $\leq 1$. However this eigenvalue can be 
complex and a corresponding eigenvector does not have to exist in the real Lie algebra $\lien$. It is however not difficult to see that also 
in this case, we can find an $x$ as claimed). 
As $\tilde{f}$ is expanding (\cite[Lemma 6]{shub69-1}), we have that $\tilde{f}^n(m)$ tends to infinity as $n$ goes to infinity for all $m\in M$ which 
are not equal to the (unique) fixed point $\tilde{m}_0$ of $\tilde{f}$. As   $\tilde{f}^n (k(x)) = k (\delta_\ast^n (x))$, this implies that 
$k(x)=\tilde{m}_0$. Moreover, the same argument applies to any point of the form $rx \in \lien$. Hence, the whole line $\R x$ is mapped 
onto the point $\tilde{m}_0$ by $k$, which 
contradicts the fact that $k$ is a proper map. So, the assumption that there exists an eigenvalue of modulus $\leq 1$ is wrong. This shows that 
$\delta_\ast$ is an expanding linear map and hence $\bar\alpha$ is an expanding affine endomorphism of the infra-nilmanifold $\Gamma\backslash N$.

\medskip

Now, since we have the information that $\bar\alpha$ is expanding, we can apply \cite[Theorem 5]{shub69-1} to conclude that $h$ is actually a homeomorphism
and hence $\bar\alpha$ and $f$ are topologically conjugate. \qed

\medskip

Note that in the above theorem it did not matter in which way we realised the fundamental group $\Gamma$ as a uniform discrete subgroup 
of $N\semi F$. It turns out that if we choose the embedding in a good way (depending on the expanding map $f$!) we can recover completely 
Gromov's result.

\begin{theorem}\label{recoverGromov}
Let $f:M\rightarrow M$ be an expanding map of a compact manifold $M$, then $f$ is topologically conjugate to 
an expanding infra-nilmanifold endomorphism.
\end{theorem} 

\underline{Proof:} We already know that $f$ is topologically conjugate to an expanding affine endomorphism $\bar\alpha$
of an infra-nilmanifold $\Gamma\backslash N$, by Theorem~\ref{main1}. So it is enough to show that this affine infra-nilmanifold endomorphism
is topologically conjugate to an expanding infra-nilmanifold endomorphism  of a possibly other infra-nilmanifold.

Let $\alpha=(d,\delta)\in \Aff(N)$ be a lift of $\bar\alpha$, hence $\alpha \Gamma \alpha^{-1}\subseteq \Gamma$. 
As $\bar\alpha$ is expanding, the map $\alpha:N\rightarrow N$ has a fixed point, say $x_0$. Now let $h:N\rightarrow N:n \mapsto x_0^{-1} n$ and 
$\Gamma'= x_0^{-1} \Gamma x_0 \subseteq \Aff(N)$. Then $\Gamma'\backslash N$ is also an infra-nilmanifold (with $\Gamma'\cong \Gamma$) and 
$h$ determines a homeomorphism $\bar{h}: \Gamma\backslash N\rightarrow \Gamma' \backslash N: \Gamma\cdot n\mapsto \Gamma' \cdot x_0^{-1}n$.
One also easily checks that $\delta \Gamma' \delta^{-1} \subseteq \Gamma'$ so that $\delta$ induces an expanding infra-nilmanifold endomorphism 
$\bar\delta$ of 
$\Gamma'\backslash N$, for which the following diagram commutes:
\[ \xymatrix{
\Gamma\backslash N\ar[r]^{\bar{h}}\ar[d]_{\bar\alpha}  &  \Gamma'\backslash N\ar[d]^{\bar\delta} \\
\Gamma\backslash N\ar[r]_{\bar{h}}  &  \Gamma'\backslash N 
}.\]
This shows that $\bar\alpha$ and hence also $f$ is topologically conjugate to the expanding infra-nilmanifold endomorphism $\delta$.
\qed

\begin{remark}
We want to stress the fact here that the infra-nilmanifold which is obtained in the theorem does not only depend on $M$, but also on the 
expanding map $f$ itself. 
\end{remark}
\section{An Anosov diffeomorphism not topologically conjugate to an infra-nilmanifold automorphism}
Analogously as in the previous section, we will show that there exists an infra--nilmanifold $M=\Gamma\backslash N$ and an Anosov diffeomorphism 
$f:M\rightarrow M$ which is not topologically conjugate to an infra--nilmanifold automorphism of $M$.

For this example, we will use a 4--dimensional
flat manifold. Again the holonomy group of the corresponding Bieberbach group will be $\Z_2$, where we 
embed $\Z_2$ as the subgroup $\{ {\mathbb I}_4, L_f \}\subseteq \GL_2(\R)$, where ${\mathbb I}_4$ is the $4\times 4$ identity matrix and 
\[ L_f= \left( \begin{array}{cccc} -1 & 0 & 0 & 0 \\ 0 & -1 & 0 &  0 \\ 0 & 0 & 1 & 0 \\ 0 & 0 & 0 & 1\end{array}\right).\]
Now, let $\Gamma$ be the torsion free, discrete  and uniform subgroup of $\R^4\semi \Z_2$ generated by 
\[ a = \left( e_1, {\mathbb I}_4 \right), 
 b = \left( e_2 , {\mathbb I}_4 \right), 
 c = \left( e_3, {\mathbb I}_4 \right), 
 d = \left( e_4, {\mathbb I}_4 \right), 
 f  = \left( \left( \begin{array}{c} 0 \\ 0 \\ \frac12 \\ \frac12 \end{array} \right), L_f \right), \]
where $e_i$ is the standard basis vector with a 1 on the $i$-th spot and 0 elsewhere. It follows that 
$\Gamma$ is a Bieberbach group, with translation subgroup $\Z^4$ generated by $a,b,c$ and $d$. 

We consider the affine map 
\[ \alpha : \R^4 \rightarrow \R^4: \left(\begin{array}{c} x \\y \\z \\ t \end{array}\right) 
\mapsto \left( \begin{array}{c}
13 x + 8 y + \frac12 \\
8 x + 5 y + \frac12\\
13 z + 8 t \\
8 z + 5 t 
\end{array} \right).\]
So 
\[ \alpha= \left(\left( \begin{array}{c}\frac12\\ \frac12 \\ 0 \\ 0    \end{array} \right), 
\left(\begin{array}{cccc}  13 & 8 & 0 & 0 \\ 8 & 5 & 0 & 0 \\ 0 & 0 & 13 & 8 \\ 0 & 0 & 8 & 5  \end{array} \right)  \right) \in \Aff(\R^4).\]
One can compute that 
\[ \alpha a \alpha^{-1} = a^{13} b^{8},\,
\alpha b \alpha^{-1} = a^8 b^5,\,
 \alpha c \alpha^{-1} = c^{13} d^{8},\,
 \alpha d \alpha^{-1} = c^{8} d^{5},\,\alpha f \alpha^{-1} = abc^{10}d^6f\]
From this, one can see that $\alpha \Gamma\alpha^{-1} = \Gamma$ and hence, $\alpha$ induces a diffeomorphism $\bar\alpha$ on $\Gamma\backslash \R^4$. Moreover, 
as the linear part of $\alpha$ only has eigenvalues of modulus different than 1, $\bar\alpha$ is an Anosov diffeomorphism.
We will show that this Anosov diffeomorphism is not topologically conjugate to an infra-nilmanifold automorphism of $\Gamma\backslash \R^4$.  Suppose on the 
contrary that $\varphi:\R^4 \rightarrow \R^4$ is a linear automorphism inducing a map $\bar\varphi$ on $\Gamma\backslash \R^4$ which is topologically conjugate 
to $\bar\alpha$. We have seen that in this case $\varphi \Gamma \varphi^{-1} = \Gamma $ and $\varphi\Z_2\varphi^{-1}=\Z_2$, which now implies that 
the matrix representation of $\varphi$ is of the form:
\[ \varphi=\left(\begin{array}{cc}  A& 0 \\ 0 & B  \end{array}\right) \mbox{ with }A,B\in \GL_2(\Z),\]
where we also used that $\varphi \Z^4 \varphi^{-1}=\Z^4$. The matrix form of $\varphi$ implies that 
\[\varphi f \varphi^{-1}= c^k d^l f \mbox{  for some }  k,l\in \Z.\] 
The fact that we suppose that $\bar\varphi$ is topologically conjugate to $\bar\alpha$ implies the existence 
of a homeomorphism $h:\Gamma\backslash \R^4\rightarrow  \Gamma\backslash \R^4$ with $\bar\alpha= h^{-1} \circ \bar\varphi \circ h$. Let 
$\bar{\alpha}_\sharp, \bar\varphi_\sharp$ and $h_\sharp$ denote the induced maps on the fundamental group $\Gamma$ of $\Gamma\backslash \R^4$.
Then, we know that, up to an inner conjugation of $\Gamma$,  $\bar{\alpha}_\sharp$ resp.~$\bar{\varphi}_\sharp$ is the same 
as conjugation with $\alpha$ resp.~$\varphi$ in $\Aff(\R^4)$ and $h_\sharp(\Z^4) =\Z^4$.  We already remark here that we will be dividing out by the derived subgroup
of $\Gamma$ in a moment, so that without any problems we can forget about the possible inner conjugations.

From $\bar\alpha= h^{-1} \circ \bar\varphi \circ h$, it follows that 
$\bar\alpha_\sharp =h^{-1}_\sharp \circ \bar\varphi_\sharp \circ h_\sharp$. We claim that this condition leads to a contradiction. To easily see this, note that 
the derived subgroup of $\Gamma$ is $[\Gamma,\Gamma]=\mbox{grp}\{a^2,b^2\}$ and the centre of $\Gamma$ is $Z(\Gamma)=\mbox{grp}\{c,d\}$.
So $Z(\Gamma)[\Gamma,\Gamma]$ is a normal subgroup of $\Gamma$ and 
\[ \Gamma/Z(\Gamma)[\Gamma,\Gamma] = \Z_2 \oplus \Z_2 \oplus \Z_2.\]
where we view the first $\Z_2$ factor as being generated by $\bar{a}$, the second factor by $\bar{b}$ and the last one by $\bar{f}$, where 
$\bar{a},\,\bar{b},\bar{f}$ denote the images of $a,\,b,\,f$ under the natural projection $\Gamma\rightarrow \Gamma/Z(\Gamma)[\Gamma,\Gamma]$.
Any automorphism of $\Gamma$ induces an automorphism of $\Gamma/Z(\Gamma)[\Gamma,\Gamma]$, which can also be seen as a linear map 
of the 3-dimensional vector space $\Z_2^3$ over the field $\Z_2$. So, we can represent the induced automorphism on $\Gamma/Z(\Gamma)[\Gamma,\Gamma]$ by 
means of a matrix in $\GL_3(\Z_2)$. 

From the conjugation relations given above, we see that $\bar\alpha_\sharp(\bar{a})=\bar{a}^{13}\bar{b}^8=\bar{a}$, $\bar\alpha_\sharp(\bar{b})=\bar{b}$ and 
$\bar\alpha_\sharp(\bar{f})= \bar a \bar b \bar f$, so the corresponding matrix in $\GL_3(\Z_2)$ is 
\[ M_\alpha = \left ( \begin{array}{ccc} 1 & 0 & 1 \\ 
0 & 1 & 1 \\
0 & 0 & 1 \end{array} \right).\]
Analogously, one can see that the matrix representations of the linear automorphisms induced by $\bar\varphi_\sharp$ and $h_\sharp$ are
of the form
\[ M_\varphi=\left( \begin{array}{ccc}
a_1 & a_2 & 0 \\
a_3 & a_4 & 0 \\
0 & 0 & 1 \end{array}\right) \mbox{ and } M_h=\left( \begin{array}{ccc}
h_1 & h_2 & h_3 \\
h_4 & h_4 & h_6 \\
0 & 0 & 1 \end{array}\right)\]
where the $a_i$ are obtained by reducing the entries of $A$ modulo 2. Now, the relation $\bar\alpha_\sharp=h^{-1}_\sharp \circ \bar\varphi_\sharp \circ h_\sharp$
implies that 
\[ M_\alpha = M_h^{-1} M_\varphi M_h.\]
By focussing on the upper left $2\times 2$ corner, one immediately gets that $M_\varphi={\mathbb I}_3$. But this then implies that also 
$M_\alpha= {\mathbb I}_3$ which is clearly a contradiction. 

\medskip

This example casts a lot of doubts on the main result of \cite{mann74-1} (Theorem C). 
Note that \cite{mann74-1} does not really contain a proof for Theorem C, but refers to the proof 
of Franks' Theorem for Anosov diffeormorphisms on tori \cite[Theorem 1]{fran69-1}. There is, up to my 
knowledge, indeed nothing wrong with \cite[Theorem 1]{fran69-1} or its proof, but to be able to generalize this to the
class of infra-nilmanifolds, it is assumed in \cite{mann74-1} (see the sentence immediately after the 
statement of Theorem A on page 423) that each homotopy class of maps from an infra-nilmanifold to itself inducing 
a hyperbolic automorphism of the fundamental group, contains a hyperbolic infra-nilmanifold automorphism.
In \cite{mann74-1}, the author refers to the wrong result of Auslander for this, but even with the use of 
Corollary~\ref{correctAusl} of the current paper, the claim does not follow. 

\medskip

In fact, the example above shows that this is not correct
and one really needs also to consider hyperbolic affine automorphisms! Of course, an affine automorphism
$\bar\alpha$ is hyperbolic if $\alpha$  (or the linear part of $\alpha$) does not have any eigenvalue 
of modulus 1.

\medskip

Unfortunately, I have not been able to give an alternative proof for the analogous version 
of \cite[Theorem C]{mann74-1} for the case of affine automorphisms. So, we are left with the 
following open question:
\begin{question}
Let $f:M\rightarrow M$ be an Anosov diffeomorphism of an infra-nilmanifold. Is it true that $f$ 
is topologically conjugate to a hyperbolic affine automorphism of the infra-nilmanifold $M$?
\end{question} 
It is very tempting to believe that the answer to this question is indeed positive. In fact,
for nilmanifolds, the arguments of A.~Manning in \cite{mann74-1} are correct (every map on 
a nilmanifold is homotopic to a nilmanifold endomorphism) and so a correct partial version of \cite[Theorem C]{mann74-1} is
\begin{theorem} 
Any Anosov diffeomorphism of a nilmanifold $M$ is topologically conjugate to a
hyperbolic nilmanifold automorphism.
\end{theorem} 
So for nilmanifolds there is no need to consider affine maps (this is also true for expanding maps).

\medskip

More generally one can even ask whether or not it is true that an Anosov diffeomorphism on 
any given compact manifold $M$ is conjugate to a hyperbolic affine automorphism of an infra-nilmanifold.
For this, it would be very useful to have a generalization of \cite[Theorem 2.1]{fran70-1} to the
case of hyperbolic affine automorphisms. However, the proof of \cite[Theorem 2.1]{fran70-1} is very dependent
on the fact that the lift of an infra-nilmanifold automorphism is really an automorphism of the 
covering Lie group and it seems rather impossible to generalize this approach to the case 
of affine automorphisms.

\medskip

Recently, there has been a lot of interest in the existence question of Anosov diffeomorphisms on 
infra-nilmanifolds (e.g\ \cite{dv08-1,gorb05-1,laur03-1,lw09-1,main06-1,mw07-1,payn09-1}). Often,
one refers to \cite[Theorem C]{mann74-1} to reduce the question to a pure algebraic question.
Luckily, in case one is only dealing with nilmanifolds, there is by the above theorem  no problem at all.
On the other hand, for infra-nilmanifolds one has to be a bit more careful. However,
for the existence question, there is not really a problem.

\begin{theorem}
Let $M$ be an infra-nilmanifold. Then the following are equivalent:
\begin{enumerate}
\item $M$ admits an Anosov diffeomorphism
\item $M$ admits a hyperbolic affine automorphism
\item $M$ admits a hyperbolic infra-nilmanifold automorphism
\end{enumerate}
\end{theorem}

\underline{Proof:}
The implications $3.\Rightarrow 2.$ and $2.\Rightarrow 1.$ are obviously true, so we only have to show 
$1. \Rightarrow 3.$ Let $M=\Gamma\backslash N$ where $N$ is a connected simply connected nilpotent Lie group
and $\Gamma$ is a uniform discrete subgroup of $N\semi F$ where $F$ is a finite subgroup of $\Aut(N)$. We 
also assume that $F$ is the holonomy group of $\Gamma$. Moreover, as is explained in \cite[section 3]{dv08-1},
we can assume that any element of $F$ restricts to an automorphism of $N_\Q$ (see \cite{dv08-1} or the proof 
of Lemma~\ref{crit1} for the meaning of $N_\Q$) 
and that $\Gamma$ is actually a subgroup of $N_\Q \semi F$ (which we called a rational realization in \cite{dv08-1}).

Assume that $f:M\rightarrow M$ is an Anosov diffeomorphism. By \cite[Theorem A]{mann74-1} $f$ induces
a hyperbolic automorphism $f_\sharp :\Pi_1(M,m_0)\cong \Gamma \rightarrow \Pi_1(M,f(m_0))\cong \Gamma$.
We recall here that for different choices of isomorphisms of $\Pi_1(M,x)$ with $\Gamma$ the induced 
map $f_\sharp:\Gamma \rightarrow \Gamma $  will change by an inner automorphism of $\Gamma$.
Anyhow, the existence of an Anosov diffeomorphism of $M$ implies the existence of a
hyperbolic automorphism $\varphi=f_\sharp$ of $\Gamma$. In the second part of the proof of Theorem A in 
\cite[page 564]{dv08-1}, we show that for some positive power $\varphi^k$ there exists a $\psi\in \Aut(N)$ such that 
$\varphi^k$ is just conjugation by $\psi \in \Aff(N)$:
\[ \forall \gamma \in \Gamma: \; \varphi^k (\gamma) = (1, \psi) \gamma (1,\psi)^{-1} .\]
As $\varphi$ is a hyperbolic, the same holds for $\varphi^k$ and hence also for $\psi$ (Lemma~\ref{zelfdeeig}).
It follows that $\Psi:N\semi F \rightarrow N\semi F:\; x \mapsto (1,\psi) x (1,\psi)^{-1}$  is an 
automorphism of $N\semi F$ with $\Psi(F)=F$ and $\Psi(\Gamma)=\Gamma$. Hence, $\Psi$ determines a
hyperbolic infra-nilmanifold automorphism of $\Gamma\backslash N$. \qed 

\medskip

Actually, the proof given above also shows the following

\begin{theorem}
An infra-nilmanifold $M$ admits an Anosov diffeomorphism if and only if $\Pi_1(M)$ admits a hyperbolic 
automorphism.
\end{theorem}

Moreover, we also showed that for any Anosov diffeomorphism $f$ on a given infra--nilmanifold $M$, 
there is some positive power $f^n$ of $f$ 
such that $f^n$ is homotopic to an infra-nilmanifold endomorphism of $M$.  Actually, this is true for any 
homeomorphism of an infra-nilmanifold. We note here that this does not hold for 
expanding maps. 

\medskip

As a conclusion of this paper we can state that in the study of selfmaps of a given infra-nilmanifold, 
which play a crucial role in the theory of expanding maps and Anosov diffeomorphisms, the class of 
infra-nilmanifold endomorphisms is just not rich enough to contain at least one map from 
each homotopy class and one really should consider the more general class of affine endomorphisms on that infra-nilmanifold. 


\end{document}